\documentclass[leqno,12pt]{amsart} %leqno is the option to put formula numbers on the left side
\usepackage{amssymb}
\theoremstyle{plain} %text of this environment is typesetted in italics

\theoremstyle{definition} %text of this environment is typesetted in roman letters

\begin{document}

\title[]{CURVATURE AND TACHIBANA NUMBERS} %title of paper and the running head option

\author[]{SERGEY E. STEPANOV} %first author's name and the running head option

%\author[S. Author]{Second Author} %second author's name and the running head option

%\dedicatory{Dedicated to Professor Xxx Yyy on his sixtieth
%birthday}

%%%%%%%%%%%%%%% footnote %%%%%%%%%%%%%%%%
%In case \subjclass[2000] command is not effective
%(or the version of amsart.cls is old), write as follows instead:
%\renewcommand{\thefootnote}{\fnsymbol{footnote}}
%\footnote[0]{2000\textit{ Mathematics Subjet Classification}.
%Primary 00; Secondary 00.}
%

%%%%%%%%%%%% Authors addresses %%%%%%%%%%%%%

%%%%%%%%%%%%%%%%%%%%%%%%%%%%%%%%%%%%%%%%%

\maketitle

{\it Department of Mathematics,
\endgraf  Financial Academy under the Government of Russian
Federation, \endgraf  125468, Moscow, Russia.
\endgraf   e-mail: {Stepanov@vtsnet.ru}}

\begin{abstract}
The purpose of this paper is to define the $r$-th Tachibana number
$t_r$ of an $n$-dimensional closed and oriented Riemannian
manifold $(M,g)$ as the dimension of the space of all conformal
Killing $r$-forms for $r = 1, 2, \dots, n - 1$ and to formulate
some properties of these numbers as an analogue to properties of
the $r$-th Betty number $b_r$  of a closed and oriented Riemannian
manifold.

\end{abstract}

\section*{0. Introduction} %delete * to number this section

{\bf 0.1.} In the present paper we define Tachibana numbers of a
compact and oriented Riemannian manifold and basing on the
definition explain well-known results, including our own. This
paper contains a new view point on well-known facts of the theory
of Riemannian manifolds. Moreover in our paper we used ideas that
were reflected in classic monograph [24].

The paper is based on the author's report given at the
International Conference "Differential Equations and Topology"
dedicated to the Centennial Anniversary of L.S. Pontryagin,
Moscow, June 17-22, 2008 (see [12]).

{\bf 0.2.} We shall concerned with four subspaces of the vector
space $\mathbf{\Omega}^r(M,\mathbf{R})$ of exterior differential
$r$-forms on an $n$-dimensional closed and oriented Riemannian
manifold $(M,g)$, namely, the subspace
$\mathbf{H}^r(M,\mathbf{R})$ of harmonic differential $r$-forms,
the subspace $\mathbf{T}^r(M,\mathbf{R})$ of {\it conformal
Killing differential $r$-forms}, the subspace
$\mathbf{K}^r(M,\mathbf{R})$ of {\it co-closed conformal Killing
$r$-forms} and the subspace $\mathbf{P}^r(M,\mathbf{R})$ of {\it
closed conformal Killing $r$-forms} for any $r=1,2,\dots, n-1$.

For a  closed and orientable  $(M,g)$, we shall call the number
$t_r = \mathrm{dim}\mathbf{T}^r(M,\mathbf{R})$ -- the {\it
Tachibana number}, $k_r=\mathrm{dim}\mathbf{K}^r(M,\mathbf{R})$ --
the {\it Killing number} and
$p_r=\mathrm{dim}\mathbf{P}^r(M,\mathbf{R})$ -- the {\it planar
number} of $(M,g)$. We shall show that the Tachibana number $t_r$
are conformal scalar invariant, while the Killing numbers $k_r$
and planar numbers $p_r$ are projective scalar invariants, and
that they are dual in sense that $t_r=t_{n-r}$ and $k_r=p_{n-r}$
for any $r=1,2,\dots, n - 1$. And, in addition, we shall deduce
the following identity $t_r=k_r+p_r$ on an $n$-dimensional closed
and orientable $(M,g)$ with the constant section curvature $C >
0$. In the last section of our paper we shall formulate the
vanishing theorem for the $r$-th Tachibana number.

\section{Curvature and Betti numbers}

In this section we summarize well known properties of harmonic
forms on closed Riemannian manifolds.

{\bf 1.1.} Let $(M,g)$  be an $n$-dimensional closed and
orientable Riemannian manifold. We shall consider the space of all
$C^\infty$-sections of  the $r$-th exterior power
$\Lambda^rM:=\Lambda^r(T^*M)$ of the cotangent bundle $T^*M$ of
$M$ as the vector space $\mathbf{\Omega}^r(M, \mathbf{R})$ of all
exterior differential $r$-forms $\omega$ on the Riemannian
manifold $(M,g)$.

 The exterior differential $r$-form $\omega$ is called a {\it harmonic
$r$-form} if $\Delta\omega= 0$ where $\Delta=d^*d+dd^*$  is the
Hodge-de Rham Laplacian acting on exterior differential forms
where $d$ is the exterior derivative and $d^*$ its formally
adjoint co-differential operator (see [1], [9]). Harmonic
$r$-forms on $(M,g)$ constitute a vector space denoted by
$\mathbf{H}^r(n, \mathbf{R})$. By the Hodge theory (see [11]), the
$r$-th Betti number $b_r$ is the dimension of $\mathbf{H}^r(n,
\mathbf{R})$ on a closed and oriented Riemannian manifold $(M,g)$.

An important property of the Hodge-de Rham Laplacian $\Delta$ is
that it commutes $*\Delta =\Delta*$ with the Hodge star operator
$*$ maps a $r$-form into an $(n - r)$-form such that
$*:\mathbf{\Omega}^r(M,\mathbf{R})\rightarrow\mathbf{\Omega}^{n-r}(M,
\mathbf{R})$ is an isomorphism and
$*^2=(-1)^{p(n-p)}\mathrm{Id}_{\Lambda^rM}$ for a fix local
orientation of a manifold $M$. In particular, if $\omega$ is a
harmonic $r$-form, then $*\;\omega$ is a harmonic $(n - r)$-form
too, i.e., $*:\mathbf{H}^r(n, \mathbf{R})\cong\mathbf{H}^{n -
r}(n, \mathbf{R})$ is an isomorphism. That implies the following
equation $b_r=b_{n-r}$ which is well known as the {\it
Poincar\'{e} duality theorem} for Betti numbers.

{\bf 1.2.} For a Riemannian manifold $(M,g)$ an arbitrary $\omega
\in \mathbf{H}^r(n, \mathbf{R})$ is a harmonic $r$-form with
respect to the conformally equivalent metric $\hat{g}=e^{2f}g$ for
any differential function $f$ on $(M,g)$ if  $n = 2r$ (see [1]).
This implies that the $r$-th Betti number $b_r$ is a {\it
conformal scalar invariant} of a Riemannian manifold $(M,g)$  of
dimension $n = 2r$.

{\bf 1.3.} The {\it rough Bochner Laplacian} is defined as
$\nabla^*\nabla$ where $\nabla$ is the Levi-Civita connection
extended to the bundle $\Lambda^rM$  of exterior differential
$r$-forms and $\nabla^*$ its formally adjoint operator (see [1]).
The rough Bochner Laplacian and the Hodge-de Rham Laplacian
related by {\it the classical Bochner-Weizenbock formula}
$\Delta\omega=\nabla^*\nabla\omega+F_r(\omega)$ where $F_r$ is a
curvature term acting on differential $r$-forms (see [1]).  If we
multiply the left- and right-hand sides of the Bochner-Weizenbock
formula by $\omega$ and integrate it on the closed and oriented
manifold $(M,g)$, then we obtain the following integral formula
$\int_Mg(\Delta\omega,\omega)dv =\int_M\|\nabla \omega\|^2dv
+\int_Mg(F_r(\omega),\omega)dv$. If $\omega$ is harmonic, then
from this integral Bochner-Weizenbock formula we have the
following (see [10])
$$
0=\int_M\|\nabla \omega\|^2dv
+\frac{1}{4}\int_M\sum_{\alpha}\lambda_\alpha\|
[\Theta_\alpha,\omega]\|^2dv \eqno(1.1)
$$
 where $\lambda_\alpha$  are the eigenvalues for the standard
symmetric {\it Riemann curvature operator}
$\overline{R}:\Lambda^2(TM)\rightarrow\Lambda^2(TM)$ given by
$g(\overline{R}(X\wedge Y),Z\wedge U)=g(R(X,Y)U,Z)$ for any vector
fields $X, Y, Z, U$ and the curvature tensor $R$ of $(M,g)$ and
$\Theta_\alpha$ are the duals of eigenvectors for the curvature
operator $\overline{R}$. We say that $(M,g)$ has {\it positive
curvature operator} if the eigenvalues of $\overline{R}$ are
positive, and we denote this by $\overline{R}>0$.

 Let $(M,g)$  has
positive curvature operator, then all $\lambda_\alpha>0$  and the
formula (1.1) implies $\nabla\omega=0$  and
$[\Theta_\alpha,\omega]=0$ for all $\alpha$. In this case $\omega$
must be zero unless $r$ is $0$ or $n$ (see also [10]). Hence,
$b_r=0$ for all $r = 1, 2, \dots, n - 1$.

 In particular, for a compact
and orientable Riemannian manifold of constant positive sectional
curvature or a compact and orientable conformally flat Riemannian
manifold with positive-definite Ricci tensor, $b_r=0$  for all  $r
= 1, 2, \dots, n - 1$.

 These are most known fundamental results of the
"Bochner technique" (see [24]).

{\bf 1.4.} Let $(M, g)$ be an $n$-dimensional closed Riemannian
manifold. For such a Riemannian manifold $(M, g)$, the space
$\mathbf{\Omega}^r(M, \mathbf{R})$ has an orthogonal (with respect
to the global scalar product) decomposition
 $$
 \mathbf{\Omega}^r(M, \mathbf{R}) = \mathrm{Ker} \;\Delta  \oplus   d (C^\infty\Lambda^{r-1}M )\oplus
 d^*( C^\infty\Lambda^{r+1}M)
 $$
where $\mathrm{Ker}\;\Delta:=\mathbf{H}^r(n, \mathbf{R})$.
Moreover, we have the orthogonal decompositions (see [8])
$$
\mathrm{Ker}\;\Delta = \mathrm{Im}\; d^*\oplus \mathrm{Ker}\; d ;
\eqno (1.2)
$$
$$
 \mathrm{ Ker}\; d^* = \mathrm{Im}\; d^*\oplus \mathrm{Ker}\;\Delta ;\quad
 \mathrm{Ker}\; d = \mathrm{Im}\; d \oplus \mathrm{Ker} \;\Delta.
 \eqno(1.3)
 $$
When $(M, g)$ is a Riemannian manifold with positive curvature
operator, we have
 $$
 \mathbf{\Omega}^r(M, \mathbf{R}) = d (C^\infty\Lambda^{r-1}M )\oplus d^* ( C^\infty\Lambda^
 {r+1}M )
 $$
for all  $r = 1, 2,\dots, n - 1$ that means
$$
\mathbf{\Omega}^r(M, \mathbf{R}) = \mathbf{D}^r(M,
\mathbf{R})\oplus \mathbf{F}^r(M, \mathbf{R}) \eqno(1.4)
$$
where $\mathbf{D}^r(M, \mathbf{R})$ and $\mathbf{F}^r(M,
\mathbf{R})$ are vector spaces of closed and co-closed exterior
differential $r$-forms on $(M, g)$ respectively.

\section{Seven vector spaces of exterior differential
$r$-forms}

In this section we define seven subspaces of the space of all
exterior differential $r$-forms on an $n$-dimensional manifold
$(M,g)$ from the stand-point of first order natural linear
differential operator on differential $r$-forms for $r = 1, 2,
\dots, n - 1$.

{\bf 2.1.} More then thirty years ago Bourguignon J. P. has
investigated (see [3]) the space $\mathrm {Diff_1}(\Lambda^rM,\Im
M)$ of natural (with respect to isometric diffeomorphisms)
differential operators of order 1 determined on vector bundle
$\Lambda^rM$  of exterior differential $r$-forms and taking their
values in the space $\Im M$ of homogeneous tensor fields on
$(M,g)$.

Bourguignon J. P. has proved the existence of {\it three basis
natural operators} of the space ${\rm Diff_1}(\Lambda^rM,\Im M)$,
but only the following two $D_1$ and $D_2$ of them were
recognized. He has shown that the first operator  $D_1$ is the
exterior differential operator $d:C^\infty\Lambda^rM\rightarrow
C^\infty\Lambda^{r+1}M$ and the second operator  $D_2$ is the
exterior co-differential operator
$d^*:C^\infty\Lambda^rM\rightarrow C^\infty \Lambda^{r-1}M$
defined by the relation (see [1])
$$d^*=(-1)^{np+n+1}*d*. \eqno(2.1)
$$

 About the third basis natural operator $D_3$, it was said that
except for case $r = 1$, this operator does not have any simple
geometric interpretation. Next, for the case $r = 1$, it was
explained that the kernel of this operator consists of
infinitesimal conformal transformations in $(M,g)$.

{\bf 2.2.} In connection with this, we have received a
specification of the Bourguignon proposition and proved (see [13])
that the basis of natural differential operators consists of three
operators of following forms:
$$
D_1=\frac{ 1}{r+1}d; \quad D_2=\frac{1}{n-r+1}g\wedge d^*; \quad
D_3=\nabla-\frac{1}{r+1}d -\frac{1}{n-r+1}g\wedge d^*
$$
where $(g\wedge d^*\omega)(X_0,X_1, \dots,
X_r)=\sum\limits_{\alpha=2}^r(-1)^\alpha
g(X_0,X_\alpha)(d^*\omega)(X_1,\dots,X_{\alpha-1},X_{\alpha+1},\dots,X_r)$
for an arbitrary $r$-form $\omega$ and vector fields
$X_1,X_2,\dots,X_r$ on $M$.

The kernel of $D_1$ consists of {\it closed exterior differential
$r$-forms}, the kernel of $D_2$ consists of {\it co-closed
exterior differential $r$-forms} and kernel of $D_3$ consists of
{\it conformal Killing $r$-forms} or, in other words, conformal
Killing tensors of order $r$ (see [5]) that constitute vector
spaces $\mathbf{D}^r(M,\mathbf{R}), \mathbf{F}^r(M,\mathbf{R})$
and $\mathbf{T^r}(M,\mathbf{R})$ respectively.

{\bf Remark.} The concept of conformal Killing tensors was
introduced by S. Tachibana about forty years ago (see [22]). He
was the first who has generalized some results of a conformal
Killing vector field (or, in other words, an infinitesimal
conformal transformation) to a skew symmetric covariant tensor of
order 2 named him the conformal Killing tensor.

 {\bf 2.3.} On a closed and
oriented Riemannian manifold $(M,g)$ the condition $\omega \in ker
\;\Delta$ is equivalent to the following condition $\omega \in
ker\; d\cap ker\; d^*$, see (1.2). Therefore the condition $\omega
\in ker\; D_1\cap ker\; D_2$ characterizes the form  $\omega$ as a
harmonic form too. Hence, the space of all harmonic $r$-forms we
can define as
$\mathbf{H}^r(M,\mathbf{R})=\mathbf{D}^r(M,\mathbf{R})\cap
\mathbf{F}^r(M,\mathbf{R})$.

The condition $\omega \in ker\; D_3 \cap ker\; D_2$ characterizes
a $r$-form $\omega$ as co-closed conformal Killing form. A
co-closed conformal Killing form is called the {\it Killing
tensor} (see, for example, [5]; [22] and [23]). The space of all
co-closed conformal Killing $r$-forms we can define as
$\mathbf{K}^r(M,\mathbf{R})=\mathbf{F}^r(M,\mathbf{R})\cap
\mathbf{T}^r(M,\mathbf{R})$.

The condition $\omega \in ker\; D_1 \cap ker\; D_2$ characterizes
the form $\omega$ as a closed conformal Killing form. Sometimes
closed conformal Killing form are also called {\it planar}. The
space of all closed conformal Killing $r$-forms we can define as
 $$
\mathbf{P}^r(M,\mathbf{R})=\mathbf{D}^r(M,\mathbf{R})\cap
\mathbf{T}^r(M,\mathbf{R}).
 $$
{\bf 2.4.} Denote the vector space of exterior differential forms
of degree $r$  by  $\mathbf{\Omega}^r(M, \mathbf{R})$  and the
vector space of parallel $r$-forms on $M$ with respect to $\nabla$
by $\mathbf{C}^r(M,\mathbf{R})$. And in this case, we have the
following diagram of inclusions (see [13]; [14]):

\begin{picture}(250,45)
  \put(45,-18){$ {\bf \mathbf{\Omega}}^r(M,\mathbf{R})$}
 \put(100,0){\vector(4,3){40}}
 \put(100,-15){\vector(3,0){40}}
 \put(100,-30){\vector(4,-3){40}}
 \put(150,25){$ {\bf D}^r(M,\mathbf{R}) $}
 \put(150,-20){$ {\bf T}^r(M,\mathbf{R}) $}
 \put(150,-60){$ {\bf F}^r(M,\mathbf{R}) $}
  \put(205,20){\vector(4,-3){40}}
 \put(205,33){\vector(1,0){40}}
 \put(205,-17){\vector(4,-3){40}}
 \put(205,-14){\vector(4,3){40}}
 \put(205,-53){\vector(4,3){40}}
 \put(205,-59){\vector(4,0){40}}
\put(250,25){$ {\bf P}^r(M,\mathbf{R})$}
 \put(250,-20){$ {\bf H}^r(M,\mathbf{R})$}
 \put(250,-60){$ {\bf K}^r(M,\mathbf{R}) $}
  \put(305,30){\vector(3,-2){40}}
 \put(305,-15){\vector(1,0){40}}
 \put(305,-60){\vector(3,2){40}}
 \put(350,-20){$ {\bf C}^r(M,\mathbf{R})$}
\end{picture}

\begin{picture}(500,70)
\end{picture}

For example, in this diagram the writing $\mathbf{F}^r(M,
\mathbf{R}) \rightarrow \mathbf{K}^r(M, \mathbf{R})$ means that
the space $\mathbf{K}^r(M, \mathbf{R})$ is subspace of the space
$\mathbf{F}^r(M, \mathbf{R})$.

\section{Tachibana numbers}

{\bf 3.1.} Let $(M,g)$  be an $n$-dimensional closed and oriented
Riemannian manifold. In the paper [13] we have fined the operator
$D^*$ formally adjoint to $D$ as follows
$$
D^*=\nabla^*-\frac{1}{r+1}d^*-\frac{1}{n-r+1}d\circ trace
$$
and constructed the rough Laplacian$^2$  $D^*_3D_3$ which has the
form
$$
D^*_3 D_3=\frac{1}{r(r+1)}\left(\nabla^*\nabla-\frac{1}{r+1}d^*
d-\frac{1}{n-r+1}d d^*\right).\eqno(3.1)
$$
Moreover, we have proved the following (see [15]):
$$
\omega\in \mathbf{T}^r(M,\mathbf{R})\Leftrightarrow \omega \in Ker
(D^*D);
$$
$$
\omega\in \mathbf{K}^r(M,\mathbf{R})\Leftrightarrow \omega \in Ker
(D^*D)\cap Ker\; d^*; \eqno(3.2)
$$
$$
\omega\in \mathbf{P}^r(M,\mathbf{R})\Leftrightarrow \omega \in Ker
(D^*D) \cap Ker\; d.
$$
The operator  $D_3$ is a typical example of so-called Stein-Weiss
operator and it was in this context already considered by T.
Branson in [4]. In particular, it was shown that $D^*_3D_3$ is
elliptic, which easily follows from computing the symbol (see also
[19]). In this case from the general theory (see [1]; [9]) we
conclude that the kernel of $D^*_3D_3$ has the finite dimension on
a closed Riemannian manifold  $(M,g)$. In this situation we see
that ${\rm dim}\mathbf{T}^r(M,\mathbf{R})= \mathrm{dim}_\mathbf{R}
{(ker\; D^*_3D_3)}<\infty$. Moreover, using this we can also
conclude the following:
$\rm{dim}\mathbf{K}^r(M,\mathbf{R})<\infty$ and
$\rm{dim}\mathbf{P}^r(M,\mathbf{R})<\infty$.

{\bf Definition.} {\it For a non-negative integer $r$ $(1\leq
r\leq n-1)$, the $r$-th Tachibana number $t_r$ of the
$n$-dimensional closed and oriented Riemannian manifold  $(M,g)$
is defined as the dimension of the space of all conformal Killing
$r$-forms on $(M,g)$.}

 And, in addition, numbers $k_r={\rm dim }\; \mathbf{K}^r(M,\mathbf{R})$
  and $p_r={\rm dim }\; \mathbf{P}^r(M,\mathbf{R})$ we shall name
   as the {\it $r$-th
Killing number} and the {\it $r$-th planar number} respectively.

{\bf 3.2.} Fix a local orientation of a manifold $M$ and consider
the action of the Hodge star operator
$*:\mathbf{\Omega}^r(M,\mathbf{R})\cong\mathbf{\Omega}^{n-r}(M,\mathbf{R})$.
The following isomorphism is obvious:
 $$
  *:\mathbf{\Omega}^r (M,\mathbf{R})\cong\mathbf{\Omega}^{n-r}(M,\mathbf{R});
  \quad*:\mathbf{C}^r (M,\mathbf{R})\cong\mathbf{C}^{n-r}(M,\mathbf{R}).
 $$
Therefore, by (2.1), we obtain the isomorphism
$$
*:\mathbf{F}^r(M,\mathbf{R})\cong\mathbf{D}^{n-r}(M,\mathbf{R}).\eqno(3.3)
$$
 From (3.3), we obtain the well-known isomorphism  $*:\mathbf{H}^r(M,\mathbf{R})\cong\mathbf{H}^{n-r}(M,\mathbf{R})$  of
vector spaces of harmonic forms. The following isomorphism (see
 [6])
 $$
 *:\mathbf{T}^r(M,\mathbf{R})\cong\mathbf{T}^{n-r}(M,\mathbf{R})\eqno(3.4)
 $$
of vector spaces of conformal Killing forms is also well known
(see also [21]). From (3.3) and (3.4) we conclude the isomorphism
(see also [13])
$$
*:\mathbf{P}^r(M,\mathbf{R})\cong\mathbf{K}^{n-r}(M,\mathbf{R}).
\eqno(3.5)
$$

From (3.4) and (3.5) we obtain that the numbers  $t_r$, $k_r$  and
$p_r$ are dual in sense that  $t_r=t_{n-r}$  and $k_r=p_{n-r}$.
These equalities are an analogue of the Poincar\'{e} duality
theorem for the Betti numbers. We have the following:

{\bf Theorem.} {\it Let $(M,g)$ be a closed and oriented
Riemannian manifold of dimension $n\geq 2$. The Tachibana numbers
$t_r$ are dual in the sense that $t_r=t_{n-r}$ and $k_r=p_{n-r}$
for $r = 1, 2, \dots , n - 1$.}

We recall that on a compact oriented Riemannian manifold $(M,g)$
with positive constant curvature there are no harmonic $r$-forms
$(1\leq r \leq n-1)$. As an application of the proposition we have
the decomposition (1.4). Using this decomposition and the
condition (3.2) Kora have concluded in [6] the following
decomposition of the space of conformal Killing $r$-forms
$$
\mathbf{T}^r(M,\mathbf{R})=\mathbf{K}^{r}(M,\mathbf{R})\oplus\mathbf{P}^r(M,\mathbf{R})\eqno(3.6)
$$
 As a consequence we have the following theorem.

{\bf Theorem.} {\it Let $(M,g)$  be an $n$-dimension $(n\geq 2)$
closed and oriented Riemannian manifold with positive constant
curvature, then the Tachibana number $t_r$ has the form
$t_r=k_r+p_r$ for any  $r = 1, 2, \dots , n - 1$. }

{\bf 3.3.} One of the most impotent properties of the equation
defining conformal Killing forms is it its conformal invariance
(see [2]; [16]), i.e., if we consider the identity conformal map
${\rm id}:(M,g)\rightarrow(M,\hat{g})$ such that
$\hat{g}:=e^{2f}g$ for any differential function $f$ on $M$ then
for any conformal Killing $r$-form $\omega$ the form
$\hat{\omega}:=e^{(r+1)f}\omega$ must be a conformal Killing
$r$-form with respect to the conformally equivalent metric
$\hat{g}:=e^{2f}g$.

On the other hand, if  we consider the identity map ${\rm
id}:(M,g)\rightarrow(M,\hat{g})$ such that preserves geodesics
then for any closed (or co-closed) conformal Killing $r$-form
$\omega$ the form $\hat{\omega}:=e^{-(r+1)f}\omega$ where
$f=(n+1)^{-1}\ln \sqrt{{{\rm det}\:g/\;{\rm det}\:\hat{g}}}$ must
be a closed (or co-closed respectively) conformal Killing $r$-form
with respect to the protectively equivalent metric $\hat{g}$ (see
[16];[20]).

We can state the following theorem.

{\bf Theorem.} {\it Let $(M,g)$ be a closed and oriented
Riemannian manifold of dimension $n\geq 2$. The Tachibana numbers
$t_r$ are conformal scalar invariants of the manifold, while the
Killing and planar numbers $k_r$ and $p_r$  are projective scalar
invariants of $(M,g)$. }

This theorem is an analogue of the statement on conformal
invariant of Betti numbers.

\section{About Tachibana numbers existence}

In this section we shall give two examples of Riemannian manifolds
with the non-vanishing Tachibana numbers.

{\bf 4.1.} Let $\mathrm{E}^{n+1}$ be a Euclidian space
$\mathrm{E}^{n+1}$ endowed with an orthogonal coordinate system
$\{x^1,x^2,\dots,x^n\}$. Then an arbitrary co-closed conformal
Killing 2-form has the following components
$\omega_{i_1i_2}=A_{ki_1i_2}x^k+B_{i_1i_2}$ where $A_{ki_1i_2}$
and $B_{i_1i_2}$ are components of a skew-symmetric constant
tensor and $k,i_1,i_2=1,2,\dots,n$ (see [23]).

Let $S^n:(x^1)^2+(x^2)^2+\dots+(x^{n+1})^2=1$ be an
$n$-dimensional unit sphere of a Euclidian space
$\mathrm{E}^{n+1}$. Tachibana has proved in [23] that $\omega
_{i_1i_2}=A_{ki_1i_2}x^k$ are components of a Killing $2$-form
defined globally on $S^n$. Now using this proposition we conclude
that the Killing number $k_2$ does not equal to zero on an
arbitrary unit hypersphere $S^n$ of an $(n+1)$-dimensional
Euclidian space $\mathrm{E}^{n+1}$. From $t_2=k_2+p_2$ we have
$t_2\geq k_2>0$ and hence $t_2\neq0$ .

In turn, we have proved in [13] the following. Let $(M,g)$ be an
$n$-dimensional locally flat pseudo-Riemannian manifold with the
rectangular coordinates  $x^1,\dots,x^n$. Then an arbitrary
co-closed conformal Killing $r$-form has the following components
$\omega_{i_1i_2\dots i_r}=A_{ki_1i_2\dots i_r}x^k+B_{i_1i_2\dots
i_r}$ where $A_{ki_1i_2\dots i_r}$ and $B_{i_1i_2\dots i_r}$ are
local components of an arbitrary skew-symmetric constant tensor
and $k,i_1,i_2,\dots,i_r=1,2,\dots,n$. In the other our paper [17]
we have presented the example $\omega_{i_1i_2\dots
i_r}=A_{ki_1i_2\dots i_r}x^k$ of a co-closed conformal Killing
$r$-form $(1\leq r \leq n-1)$ defined globally on a unit
hypersphere $S^n$. Therefore, we can conclude that the Killing
number $k_r$ $(1\leq r \leq n-1)$ does not equal to zero on an
arbitrary hypersphere $S^n$ of an $(n+1)$-dimensional Euclidian
space $\mathrm{E}^{n+1}$. Consequently, we have $t_r\neq0$ because
$t_r=k_r+p_r>0$ for all  $r = 1, \dots , n - 1$.

{\bf Theorem.} {\it The Tachibana number $t_r$ $(1\leq r\leq n-1)$
does not equal to zero on an $n$-dimensional unit sphere $S^n$ of
a Euclidian space $\mathrm{E}^{n+1}$.}

{\bf 4.2.} Let now $(M_1,g_1)$  and $(M_2,g_2)$ be Riemannian
manifolds where ${\rm dim} M_{1}=1$ and ${\rm dim} M_{n-1}=n-1$.
Let $f$ be a positive smooth function on $M_1$. The {\it warped
product} $M_1\times_f M_2$ of $(M_1,g_1)$ and $(M_2,g _2)$ is the
topological product manifold $M=M_1\times M_2$ with the metric
$g=g_1\times_f g_2$ defined by $g_1\times_f g_2=\pi_1^*g_1+(f\circ
\pi_1)\pi_2^*g_2$, where $\pi_1:M_1\times M_2\rightarrow M_1$ and
$\pi_2:M_1\times M_2\rightarrow M_2$ are natural projections on
$M_1$ and $M_2$ respectively. If $(M_1,g_1)$ and $(M_2,g_2)$ are
closed and oriented Riemannian manifolds then the warped product
$M_1\times_f M_2$ is a closed and oriented Riemannian manifold
also.

For any nonvanishing differential function $\lambda:M_1\rightarrow
\mathbf{R}$, the $1$-form $\omega={\rm grad \; \lambda}$ is a
closed conformal Killing one-form or, in other words, {\it
concircular} vector field is globally defined on
$S_1\times_fS^{n-1}$ (see [7]). Therefore, we can conclude the
following: $p_1=k_{n-1}\geq 1$. Thus, we see that the Tachibana
numbers are $t_1\geq1$ and $t_{n-1}\geq1$ because $t_1\geq p_1$
and $t_{n-1}\geq k_{n-1}$ respectively.

{\bf Theorem.} {\it The Tachibana numbers $t_r$ have the following
relations $t_1\geq 1$ and  $t_{n - 1}\geq1$   $t_1$ and $t_{n -1}$
do not equal to zero on warped product $S^1\times_f S^{n -1}$.}

\section{The vanishing theorem for Tachibana numbers}

In this section, we shall prove the theorem about vanishing of the
Tachibana numbers on a closed and oriented Riemannian manifold of
negative curvature operator; the theorem will be proved with help
of  the "Bochner technique" (see [24]; [10]).

{\bf 5.1.} If we multiply the left- and right-hand sides of
relation
$$
D^*_3
D_3\omega=\frac{1}{r(r+1)}\left(\nabla^*\nabla\omega-\frac{1}{r+1}d^*
d\omega-\frac{1}{n-r+1}d d^*\omega\right)
$$
 by  $\omega$  and integrate it on the closed and
oriented manifold $(M,g)$, we obtain the following integral
formula
$$
r(r+1)\int_M\|D_3\omega\|^2 dv=\int_M\|\nabla\omega\|^2
dv-\frac{1}{r+1}\int_M\|d\omega\|^2dv-\frac{1}{n-r+1}\int_M\|d^*\omega\|^2dv.
\eqno(5.1)
$$

 On the other hand, integrating the classical
Bochner-Weizenbock formula
$\Delta\omega=\nabla^*\nabla\omega+F_r(\omega)$ on the closed and
oriented Riemannian  manifold $(M,g)$, we obtain the second
integral formula (see [10])
$$
\int_M\|d\omega\|^2 dv+\int_M\|d^*\omega\|^2
dv=\int_M\|\nabla\omega\|^2dv+\frac{1}{4}\int_M\sum\limits
_{\alpha}\lambda_\alpha\|[\Theta_\alpha,\omega]\|^2dv. \eqno(5.2)
$$
 The equation (5.1) by means of (5.2) becomes
 $$
 r(r+1)\int_M\|D_3\omega\|^2 dv=
 $$
$$
=\frac{r}{r+1}\int_M\|d\omega\|^2
dv+\frac{n-r}{n-r+1}\int_M\|d^*\omega\|^2dv-\frac{1}{4}\int_M\sum\limits
_\alpha\lambda_\alpha\|[\Theta_\alpha,\omega]\|^2dv. \eqno(5.3)
 $$

{\bf 5.2.} We say that $(M,g)$  has {\it negative curvature
operator} if the eigenvalues of $\overline{R}$  are negative, and
we denote this by $\overline{R}  < 0$. Let $(M,g)$ be a Riemannian
manifold with negative curvature operator. This means that the
inequality $\lambda_\alpha<0$ holds for any $\alpha$. Moreover, if
we suppose that $\omega$ is a conformal Killing $r$-form then
(5.3) gives us $\nabla\omega=0$ and $[\Theta_\alpha,\omega]=0$ for
all $\alpha$. In this case the $r$-form $\omega$ must be zero
unless $r$ is $0$ or $n$ (see [10]). It means, that there exists
no conformal Killing $r$-form other then zero form on an
$n$-dimensional compact and oriented Riemannian manifold with
negative curvature operator (see also [18]). Hence, $t_r$ must be
zero unless $r$ is $0$ or $n$. Then we obtain the following
theorem.

{\bf Theorem.} {\it The Tachibana number $t_r$ $(1\leq r \leq
n-1)$ is equal to zero on an $n$-dimensional closed and oriented
Riemannian manifold $(M,g)$ with negative curvature operator.}

Using this theorem we can conclude that for a compact and
orientable Riemannian manifold of constant negative sectional
curvature or a compact and orientable conformally flat Riemannian
manifold with negative-definite Ricci tensor, $t_r=0$  $(r = 1,
2,\dots, n - 1)$.

 This theorem is an analogue of  the vanishing
theorem for Betti numbers.

\section{Two open problems}

{\bf 6.1.} To find the Tachibana numbers $t_r$  for $r = 1,
2,\dots , n - 1$ on an $n$-dimensional unit sphere $S^n$ of a
Euclidian space $\mathrm{E}^n+1$.

{\bf 6.2.} Are true or false the following equalities
$t_r=k_r+p_r$ for $r = 1, 2, \dots , n - 1$ on an $n$-dimensional
closed and oriented Riemannian manifold $(M,g)$ with positive
curvature operator whose metric $g$ has not constant positive
sectional curvature?

%%%%%%%%%%%% References %%%%%%%%%%%%%
%%
%<Author name> is written as Initial of Given Name, and Family Name.
%<Title> is written in roman letters.
%<Journal name> should be abbreviated according to
% the MR Serials Abbreviations List of Mathematical Reviews:
% (Abbreviations of Names of Serials; http://www.ams.org/mr-database)
%For <Pages>, use en-dash "--" between page numbers.
%%

\end{document}